\documentclass[11pt,leqno]{article}
\usepackage{amsmath,amstext}

\date{25 July 1998}

\parfillskip=\parindent plus1fil 

\renewcommand{\phi}{\varphi}

\newcommand{\ip}[1]{\langle #1 \rangle}

\newcommand{\wt}{\widetilde}    
\newcommand{\tl}{\tilde}
\newcommand{\til}{\tl{\ }}

\newcommand{\bibref}[1]{[\ref{#1}]}

\renewcommand{\emph}{\textsl}

\newcommand{\proof}{\textsc{Proof:} \hspace{.4em minus .2em}}

\newcommand{\qed}{\parfillskip=0pt plus1fil\nolinebreak\hfill$\rule{1ex}{1ex}$%
\par\addvspace{12pt plus 3pt minus 3pt}\parfillskip=\parindent plus1fil}

\title{\sffamily\bfseries\Large Compact Operators via the Berezin Transform}
\author{\sffamily\bfseries Sheldon Axler \and \sffamily\bfseries Dechao Zheng}

\newenvironment{keeptogether}{\pagebreak[0]\samepage}{}

\numberwithin{equation}{section}

\newtheorem{Theorem}[equation]{\sffamily\bfseries Theorem}
\newtheorem{Proposition}[equation]{\sffamily\bfseries Proposition}
\newtheorem{Lemma}[equation]{\sffamily\bfseries Lemma}
\newtheorem{Corollary}[equation]{\sffamily\bfseries Corollary}
\newtheorem{Conjecture}[equation]{\sffamily\bfseries Conjecture}

\newenvironment{theorem}{\begin{keeptogether}\begin{Theorem}\slshape}%
  {\end{Theorem}\end{keeptogether}}
  {\end{Proposition}\end{keeptogether}}
\newenvironment{lemma}{\begin{keeptogether}\begin{Lemma}\slshape}%
  {\end{Lemma}\end{keeptogether}}
\newenvironment{corollary}{\begin{keeptogether}\begin{Corollary}\slshape}%
  {\end{Corollary}\end{keeptogether}}
  {\end{Conjecture}\end{keeptogether}}

\newcounter{referencec}

\newcounter{subtheoremc}
\newcounter{subsubtheoremc}

\newenvironment{subtheorem}{\begin{list}{{\normalfont
(\alph{subtheoremc})}\hfill}{\usecounter{subtheoremc}
\setlength{\labelwidth}{.4in}
\setlength{\labelsep}{0pt}
\setlength{\leftmargin}{.4in}
}}{\end{list}}

\newenvironment{subsubtheorem}{\begin{list}{{\normalfont
(\roman{subsubtheoremc})}\hfill}{\usecounter{subsubtheoremc}
\setlength{\labelwidth}{.4in}
\setlength{\labelsep}{0pt}
\setlength{\leftmargin}{.4in}
}}{\end{list}}

\makeatletter

\renewcommand{\@startsection}[6]{\if@noskipsec \leavevmode
\fi
   \par \@tempskipa #4\relax
   \@afterindenttrue
   \ifdim \@tempskipa <\z@ \@tempskipa -\@tempskipa \@afterindenttrue\fi
   \if@nobreak \everypar{}\else
     \addpenalty{\@secpenalty}\addvspace{\@tempskipa}\fi \@ifstar
     {\@ssect{#3}{#4}{#5}{#6}}{\@dblarg{\@sect{#1}{#2}{#3}{#4}{#5}{#6}}}}

\renewcommand{\section}{\@startsection {section}{1}{\z@}%
                                   {-3.5ex \@plus -1ex \@minus -.2ex}%
                                   {2.3ex \@plus.2ex}%
                                   {\reset@font\large\sffamily\bfseries}}

\makeatother

\renewcommand{\[}{\parfillskip=0pt plus1fil$$}
\renewcommand{\]}{$$\parfillskip=\parindent plus1fil}

\makeatletter
\renewenvironment{equation}{\parfillskip=0pt plus1fil%
\refstepcounter{equation}%
$$%
}{%
\leqno\tagform@{\theequation}%
$$%
\parfillskip=\parindent plus1fil}%

\newcommand{\simpletag}[1]{\def\@currentlabel{#1}\def\theequation{#1}%
  \addtocounter{equation}{-1}}
\makeatother

\begin{document}
\maketitle
\thispagestyle{empty}

\begin{abstract}
\noindent
In this paper we prove that if $S$ equals a finite sum of finite products of Toeplitz
operators on the Bergman space of the unit disk, then $S$ is compact if and only if
the Berezin transform of $S$ equals $0$ on $\partial D$.  This
result is new even when
$S$ equals a single Toeplitz operator.  Our main result can be used to prove,
via a unified approach, several previously known results about compact Toeplitz
operators, compact Hankel operators, and appropriate products of these operators.
\end{abstract}

\section{Introduction}

{}Let $dA$ denote Lebesgue area measure on the unit disk $D$, normalized so that the
measure of~$D$ equals $1$. The \emph{Bergman space} $L^{2}_{a}$ is the Hilbert
space consisting of the analytic functions on~$D$ that are also in~$L^{2}(D, dA)$.
For $z \in D$, the \emph{Bergman reproducing kernel} is the function $K_z \in L^2_a$
such that\footnote{Both authors
were partially supported by the National Science Foundation.}
\[
f(z) = \ip{f, K_z}
\]
for every $f \in L^2_a$.  The \emph{normalized Bergman reproducing kernel} $k_z$ is the
function $K_z/\|K_z\|_2$.  Here, as elsewhere in this paper, the norm $\|\ \|_2$ and
the inner product $\ip{\ \,,\ }$ are taken in the space $L^2(D, dA)$.

For $S$ a bounded operator on $L^2_a$, the \emph{Berezin transform} of $S$ is the
function $\tl{S}$ on $D$ defined by 
\[
\tl{S}(z)=\ip{Sk_{z}, k_{z}}.
\]

For $u\in L^\infty(D, dA)$, the \emph{Toeplitz operator} $T_u$ with symbol $u$ is the
operator on $L^2_a$ defined by $T_{u}f = P(uf)$, where $P$ is the orthogonal projection
from $L^{2}(D, dA)$ onto $L_{a}^{2}$.

In this paper we prove that if $S$ equals a finite sum of finite products of Toeplitz
operators, then $S$ is compact if and only if $\tl{S}(z) \to 0$ as $z \to \partial D$. 
This result is new even when $S$ equals a single Toeplitz operator~$T_u$.  Our
main result can be used to prove, via a unified approach, several previously known
results about compact Toeplitz operators, compact Hankel operators, and appropriate
products of these operators.

A common intuition is that for operators on the Bergman space ``closely associated
with function theory'', compactness is equivalent to having vanishing Berezin transform
on~$\partial D$.  Our main result shows that this intuition is correct if 
``closely associated with function theory'' is interpreted to mean that the operator
is a finite sum of finite products of Toeplitz operators.

Section~2 of this paper contains a precise statement of our theorem along with
a discussion of some consequences and examples.  Section~3 contains
three lemmas that will be used in the proof of the theorem.  Section~4 contains the
proof of the theorem.

\section{Discussion of the Theorem}

A nice survey of previously known results connecting the Berezin transform with
Toeplitz operators (on both the Hardy space and the Bergman space) can be found in
\bibref{St2}.

Before stating our theorem, we need to introduce some notation.  For $z \in D$, let
$\phi_{z}$ be the analytic map of $D$ onto $D$ defined by
\begin{equation} \label{82}
\phi_{z}(w)=\frac{z-w}{1-\bar{z}w}.
\end{equation}%
A simple computation shows that $\phi_z \circ \phi_z$ is the identity function on~$D$.

For $z \in D$, let $U_z \colon L^2_a \to L^2_a$ be the unitary operator defined by
\[
U_z f = (f \circ \phi_z) {\phi_z}'.
\[
Notice that ${U_z}^* = {U_z}^{-1} = U_z$, so $U_z$ is actually a self-adjoint unitary
operator.

For $S$ a bounded operator on $L^2_a$, define $S_z$ to be the bounded operator on
$L^2_a$ given by conjugation with $U_z$:
\[
S_z = U_z S U_z.
\]

Although we are concerned only with operators on the Bergman space $L^2_a$, we will
need to make use of other norms.  For a measurable function $u$ on
$D$ and $1 \le p \le \infty$, let $\|u\|_p$ denote the usual norm of $u$ in
$L^p(D, dA)$.

Now we are ready to state our main result.  Although we will discuss some
consequences of this theorem in this section, we will not complete the proof of
this theorem until Section~4.

\begin{theorem} \label{37}
Suppose\/ $S$ is a finite sum of operators of the form\/ $T_{u_1} \dots T_{u_n}$, where
each\/ $u_j \in L^\infty(D,dA)$.  Then the following are equivalent:
\begin{subtheorem}
\item
$S$ is compact;

\item
$\|Sk_z\|_2 \to 0$ as\/ $z \to \partial D$;

\item
$\tl{S}(z) \to 0$ as\/ $z \to \partial D$;

\item
$S_z 1 \to 0$ weakly in\/ $L^2_a$ as\/ $z \to \partial D$;

\item
$\|S_z1\|_2 \to 0$ as\/ $z \to \partial D$;

\item
$\|S_z 1\|_p \to 0$ as\/ $z \to \partial D$ for
every\/ $p \in (1,\infty)$.
\end{subtheorem}
\end{theorem}

Our main concern is the equivalence of conditions (a) and (c) above, but the other
conditions are also of interest.  Furthermore, the other conditions are needed as
intermediate steps in the proof that (c) implies~(a).  One curious feature of the
proof is that even though we are only interested in the Bergman space $L^2_a$, we
will need to use condition (f) with $p=6$ as part of the chain of implications from
(c) to~(a), as the reader will note when we present the proof in Section~4.

The Berezin transform  $\tl{u}$ of a function $u \in L^\infty(D, dA)$ is defined to
be the Berezin transform of the Toeplitz operator $T_u$.  In other words,
$\tl{u} = \wt{T_u}$.  Note that $\tl{u}(z) = \wt{T_u}(z) = \ip{T_uk_z, k_z} =
\ip{P(uk_z), k_z} = \ip{uk_z, k_z}$ for each $z \in D$.  Expressing this as an
integral, we obtain the formula
\begin{equation} \label{81}
\tl{u}(z) = \int_D u(w) |k_z(w)|^2\,dA(w)
\end{equation}%
for each $z \in D$.  Because $\|k_z\|_2 = 1$, this formula shows that $\tl{u}(z)$ is a
weighted average of~$u$.

As is well known, the explicit formulas for the
reproducing kernel and the normalized reproducing kernel are given by
\begin{equation} \label{84}
K_z(w) = \frac{1}{(1 - \bar{z}w)^2}\ , \quad k_z(w) = \frac{1 - |z|^2}{(1 -
\bar{z}w)^2}
\end{equation}%
for $z, w \in D$.  Along with the formula \eqref{82} for $\phi_z$, this shows that
$|{\phi_z}'(w)| = |k_z(w)|$ for all $z, w \in D$.  Thus making the change of variables
$\lambda = \phi_z(w)$ in~\eqref{81}, we have $dA(\lambda) = |k_z(w)|^2\,dA(w)$ and $w =
\phi_z(\lambda)$ (because $\phi_z$ is its own inverse under composition), transforming
\eqref{81} into the formula
\[
\tl{u}(z) = \int_D (u \circ \phi_z)(\lambda)\,dA(\lambda)
\]
for each $z \in D$.

The corollary below is an immediate consequence of the equivalence of conditions (a)
and (c) in Theorem~\ref{37}.  Previously, the best results known along these lines
were special cases that required additional hypotheses on~$u$.  Specifically, Zhu
(\bibref{Zhu}, Theorem~B) proved the result below under the additional hypothesis that
$u$ is a nonnegative function.  More recently, Korenblum and Zhu \bibref{KZh} proved
the result below under the additional hypothesis that $u$ is a radial function (meaning
$u(z) = u(|z|)$ for all $z \in D$) and Stroethoff \bibref{Str} proved it under the
additional hypothesis that $u$ is uniformly continuous with respect to the hyperbolic
metric.

\begin{corollary} \label{83}
If\/ $u \in L^\infty(D, dA)$, then\/ $T_u$ is compact if and only if\/
$\tl{u}(z) \to 0$ as $z \to \partial D$.
\end{corollary}

Now we turn to a discussion of how Theorem~\ref{37} can be used prove several known
results about Toeplitz and Hankel operators.  The advantage of using Theorem~\ref{37}
is that it allows a single technique to be used in several different contexts,
replacing more specialized techniques and estimates that might work only in one
context.

We begin with a characterization of the compact Toeplitz operators.  Actually
Corollary~\ref{83} above provides the most useful characterization of the compact
Toeplitz operators, in the sense that the condition $\tl{u}(z) \to 0$ as $z \to
\partial D$ is easier to verify in practice than the conditions we will discuss
below.  However, the other conditions below are also useful, and we want to
demonstrate how Theorem~\ref{37} gives proofs of these results.

Suppose $u \in L^\infty(D, dA)$.  Zheng proved (\bibref{Zhe}, Theorem~4) that the
following are equivalent:
\begin{subsubtheorem}
\item
$T_u$ is compact;

\item
$\|T_u k_z\|_2 \to 0$ as $z \to \partial D$;

\item
$\|P(u \circ \phi_z)\|_2 \to 0$ as $z \to \partial D$;

\item
$\|P(u \circ \phi_z)\|_p \to 0$ as $z \to \partial D$ for every $p \in (1, \infty)$.
\end{subsubtheorem}
To prove this using Theorem~\ref{37}, let $S = T_u$.  The
equivalence of (i) and (ii) above is a special case of the
equivalence of conditions (a) and (b) in Theorem~\ref{37}.  An easy calculation (see
Lemma~8 of \bibref{AxC}) shows that
$S_z = T_{u \circ \phi_z}$, so $S_z 1 = P(u \circ \phi_z)$.  Thus the equivalence of
(i), (iii), and (iv) above follows from the equivalence of conditions (a), (e), (f) in
Theorem~\ref{37}, completing our proof of Zheng's result.

For the next application of Theorem~\ref{37}, we need a formula for $S_z$ when $S$ is a
finite product of Toeplitz operators.  If
$u_1, \dots, u_n \in L^\infty(D, dA)$, then
\begin{equation} \label{26}
U_z T_{u_1} \dots T_{u_n} U_z =  T_{u_1 \circ \phi_z} \dots T_{u_n \circ \phi_z}
\end{equation}%
because we can write the operator on the left side as
\[
(U_z T_{u_1} U_z)(U_z T_{u_2} U_z) \dots (U_z T_{u_n} U_z)
\]
and then use the formula $U_z T_u U_z = T_{u \circ \phi_z}$ (see Lemma~8 of
\bibref{AxC}).

Now we show how the compact Hankel operators can be characterized by using
Theorem~\ref{37}.  Again suppose that $u \in L^\infty(D, dA)$.  The
\emph{Hankel operator} with symbol
$u$ is the operator $H_u$ from $L^2_a$ to $L^2(D,dA) \ominus L^2_a$ defined by
$H_u f = (1 - P)(uf)$. 

Stroethoff and Zheng proved (\bibref{St2},
Theorem~6, and \bibref{Zhe}, Theorem~3) that the following are equivalent:
\begin{subsubtheorem}
\item
$H_u$ is compact;

\item
$\|H_u k_z\|_2 \to 0$ as $z \to \partial D$;

\item
$\|u \circ \phi_z - P(u \circ \phi_z)\|_2 \to 0$ as $z \to \partial D$.
\end{subsubtheorem}
To prove this using Theorem~\ref{37}, let $S = {H_u}^*H_u$.  We need to work
with ${H_u}^*H_u$ instead of $H_u$ because $H_u$ is not a finite sum of finite products
of Toeplitz operators.  However, the definitions of Hankel and Toeplitz operators
easily lead to the identity
\begin{equation} \label{74}
{H_u}^*H_u = T_{|u|^2} - T_{\bar{u}} T_u,
\end{equation}%
so Theorem~\ref{37} can be applied to~$S$.  The equivalence of conditions (a) and
(c) in Theorem~\ref{37} shows that ${H_u}^*H_u$ is compact (which is equivalent to
$H_u$ being compact) if and only if
$({H_u}^*H_u)\til(z) \to 0$ as $z \to \partial D$.  However
\[
({H_u}^*H_u)\til(z) = \ip{{H_u}^*H_u k_z, k_z} = {\|H_u k_z\|_2}^2,
\]
so we conclude that $H_u$ is compact if and only if $\|H_u k_z\|_2 \to 0$ as
\mbox{$z \to \partial D$}.  In other words, conditions (i) and (ii) above are
equivalent.  To get to condition~(iii), note that from \eqref{26} and~\eqref{74}, we
have $S_z = {H_{u \circ \phi_z}}^*H_{u \circ \phi_z}$.  The equivalence of conditions
(a) and (b) in Theorem~\ref{37} thus shows that $H_u$ is compact if and only if
$\|{H_{u \circ \phi_z}}^*H_{u \circ \phi_z} 1\|_2 \to 0$ as $z \to \partial D$. 
Because
\[
{\|H_{u \circ \phi_z}1\|_2}^2 \le \|{H_{u \circ \phi_z}}^*H_{u \circ \phi_z}1\|_2
\le \|{H_{u \circ \phi_z}}^*\| \|H_{u \circ \phi_z}1\|_2,
\]
we see that $H_u$ is compact if and only if $\|H_{u \circ \phi_z}1\|_2 \to 0$ as
$z \to \partial D$.  Finally, note that
$H_{u \circ \phi_z}1 = u \circ \phi_z - P(u \circ \phi_z)$.  Thus conditions (i) and
(iii) above are equivalent, completing our proof of Stroethoff's and Zheng's result.

As a special case of the last result, suppose $u = \bar{f}$, where $f$ is a bounded
analytic function on~$D$.  Then $P(\bar{f} \circ \phi_z)$ equals the constant
function $\bar{f}(z)$.  Thus from the equivalence of conditions (i) and (iii) above we
conclude that $H_{\bar{f}}$ is compact if and only if
$\|f \circ \phi_z - f(z)\|_2 \to 0$ as $z \to \partial D$.  This last condition holds
if and only if $f$ is in the little Bloch space (see Theorem~2 of \bibref{Axl}). 
Thus we have recovered Axler's result (\bibref{Axl}, Theorem~7) that $H_{\bar{f}}$ is
compact if and only if $f$ is in the little Bloch space.

We emphasize that the proofs just given of the results due to Zheng, Stroethoff and
Zheng, and Axler should not be regarded as entirely new proofs of these results.  Some
techniques from the original proofs (along with some techniques used by Stroethoff and
Zheng in~\bibref{StZ}) have been incorporated into part of our proof of
Theorem~\ref{37} (specifically, some of these techniques are used in our proof
that condition~(f) of Theorem~\ref{37} implies
condition~(a) of Theorem~\ref{37}).  The techniques just mentioned could be used to
prove that conditions (a) and (b) in Theorem~\ref{37} are equivalent, but they do not
seem powerful enough to prove that conditions (a) and (c) are equivalent (this is our
main goal).  Thus the proof of Theorem~\ref{37} requires some additional techniques,
which we introduce in Section~4.

Our paper~\bibref{AxZ} gives two additional applications of
Theorem~\ref{37}, including a characterization of the bounded analytic functions
$f, g$ such that \mbox{$T_f T_{\bar g} - T_{\bar g} T_f$} is compact.

Now we turn to a discussion of the role of the hypothesis of Theorem~\ref{37}, which
requires that $S$ be a finite sum of finite products of Toeplitz operators.  Two of
the trivial implications in Theorem~\ref{37}, namely (a) $\Rightarrow$ (b) and (b)
$\Rightarrow$ (c), obviously hold for arbitrary bounded operators $S$ on~$L^2_a$.  An
examination of the proof (see Section~4) shows that one of the difficult
implications, namely (c) $\Rightarrow$ (d), also holds for arbitrary bounded operators
$S$ on~$L^2_a$.  The proofs of the other implications all use the hypothesis that $S$
is a finite sum of finite products of Toeplitz operators.

The paragraph above suggests that we should consider whether
Theorem~\ref{37} holds for arbitrary bounded operators $S$ on~$L^2_a$.  Are the key
conditions (a), (b), and (c) of Theorem~\ref{37} equivalent for arbitrary bounded
operators $S$ on~$L^2_a$?  Unfortunately this question has a negative answer.  In
fact, we will now give some examples to show that no two of the conditions (a), (b),
and (c) of Theorem~\ref{37} are equivalent for arbitrary bounded operators $S$
on~$L^2_a$.  Thus the hypothesis that $S$ is a finite sum of finite products of
Toeplitz operators (or some appropriate substitute) is needed in Theorem~\ref{37}.

For our examples we will need a power series formula for the Berezin transform of a
bounded operator $S$ on~$L^2_a$.  From
\eqref{84} we get
\[
k_z(w) = (1 - |z|^2) \sum_{m=0}^\infty (m+1) \bar{z}^m w^m
\]
for $z, w \in D$.  To compute $\tl{S}(z)$, which equals $\ip{Sk_z, k_z}$, first
compute $Sk_z$ by applying
$S$ to both sides of the equation above, and then take the inner product with $k_z$,
again using the equation above, to obtain
\begin{equation} \label{70}
\tl{S}(z) = 
(1 - |z|^2)^2
\sum_{m,n=0}^\infty (m+1) (n+1) \ip{ Sw^m, w^n } \bar{z}^m z^n.
\end{equation}%

Consider the operator $S$ on $L^2_a$ defined by
\[
S\big( \sum_{n=0}^\infty a_n w^n \big) =  \sum_{n=0}^\infty a_{{}_{2^n}} w^{2^n}.
\]
Clearly $S$ is a self-adjoint projection with infinite-dimensional range.  Thus $S$
is not compact.  Furthermore,
\begin{align*}
{\|Sk_z\|_2}^2 &= \ip{Sk_z, k_z} \\
&= \tl{S}(z) \\
&= (1 - |z|^2)^2 \sum_{n=0}^\infty (2^n + 1) (|z|^2)^{2^n},
\end{align*}
where the first equality holds because $S$ is a self-adjoint projection and the last
equality follows from~\eqref{70}.  Because
$(1-t)^2 \sum_{n=0}^\infty (2^n + 1) t^{2^n} \to 0$ as $t \to 1^-$ (we leave the
verification of this limit as an exercise), the equation above shows that
${\|Sk_z\|_2}^2 = \tl{S}(z) \to 0$ as $z \to \partial D$.  This shows that neither
condition~(b) nor condition~(c) of Theorem~\ref{37} implies condition~(a) for
arbitrary operators on~$L^2_a$.

Now consider the operator $S$ on $L^2_a$ defined by
\[
S\big( \sum_{n=0}^\infty a_n w^n \big) =  \sum_{n=0}^\infty (-1)^n a_n w^n.
\]
Clearly $S$ is a unitary operator, so $\|Sk_z\|_2 = 1$ for all $z \in D$. 
Thus $\|Sk_z\|_2 \not\to 0$ as $z \to \partial D$.
However, a calculation (use \eqref{70}) shows that
\[
\tl{S}(z) = \frac{(1 - |z|^2)^2}{(1 + |z|^2)^2}
\]
for $z \in D$.  Thus $\tl{S}(z) \to 0$ as $z \to \partial D$.  Hence $S$ satisfies
condition~(c) of Theorem~\ref{37} but not condition~(b), showing that
condition~(c) does not imply condition~(b) for arbitrary operators on~$L^2_a$.  This
example is based on a similar example provided by Nordgren and Rosenthal
\bibref{NoR} in the context of Hardy spaces.

After this paper appeared in preprint form, Miroslav Engli\u{s} (private
communication) showed that a conjectured quantitative improvement of Theorem~\ref{37}
is false.  Specifically, he shows that the essential norm of $T_u$ and
$\limsup_{z \to \partial D} |\tl{u}(z)|$ are not equivalent for
$u \in L^\infty(D, dA)$.

\section{Some Lemmas}

This section contains three lemmas that will be used in the proof of
Theorem~\ref{37}.  We begin by stating a simple lemma describing how the Berezin
transform acts under composition with the maps $\phi_z$.  Recall that $K_z$ denotes the
reproducing kernel on $L^2_a$ and that $k_z$ denotes the normalized reproducing
kernel $K_z/\|K_z\|_2$.

\begin{lemma} \label{77}
If\/ $S$ is a bounded operator on\/ $L^2_a$ and\/ $z \in D$, then\/ $\tl{S}
\circ
\phi_z =
\wt{S_z}$.
\end{lemma}

The proof of the lemma above is a calculation that we leave to the reader as an
exercise.

Note that $(S_w)^* = (S^*)_w$ for every bounded operator $S$ on $L^2_a$
and every $w \in D$.  Thus the expression $S_w^*$, which appears in \eqref{51} below,
can be interpreted either way.

For $S$ a bounded operator on $L^2_a$ and $z, w \in D$, the function
$S_z 1$ or the function $S_w^*1$ may not be in $L^6(D, dA)$.  In that case, the right
side of \eqref{50} or \eqref{51} below will be infinity, making the corresponding
inequality trivially true.

The appearance of the $L^6$ norm in the next lemma may seem
strange.  However, the lemma becomes false if the $6$ is replaced by~$2$, because
otherwise we would be able to prove that conditions (a) and (b) in Theorem~\ref{37}
are equivalent for arbitrary bounded operators $S$ on~$L^2_a$ (which we have already
shown is false).

\begin{lemma} \label{58}
There is a constant\/ $c < \infty$ such that if\/ $S$ is a bounded operator
on\/ $L^2_a$, then
\begin{equation} \label{50}
\int_D \frac{|(SK_z)(w)|}{\sqrt{1 - |w|^2}}\,dA(w) \le
\frac{c \|S_z 1\|_6}{\sqrt{1 - |z|^2}}
\end{equation}%
for all\/ $z \in D$ and
\begin{equation} \label{51}
\int_D \frac{|(SK_z)(w)|}{\sqrt{1 - |z|^2}}\,dA(z) \le
\frac{c \|S_w^*1\|_6}{\sqrt{1 - |w|^2}}
\end{equation}%
for all\/ $w \in D$.
\end{lemma}

\proof
An easy calculation shows that
\begin{equation} \label{86}
U_z 1 = (|z|^2 - 1) K_z.
\end{equation}%

To prove \eqref{50}, fix a bounded operator $S$ on $L^2_a$ and fix $z \in D$.  We have
\[
SK_z = \frac{SU_z 1}{|z|^2-1} = \frac{U_z S_z 1}{|z|^2-1}
= \frac{((S_z 1)\circ \phi_z){\phi_z}'}{|z|^2-1},
\]
where the first equality comes from \eqref{86}, the second equality comes from the
definition of~$S_z$, and the third equality comes from the definition of~$U_z$.  Thus
\[
\int_D \frac{|(SK_z)(w)|}{\sqrt{1 - |w|^2}}\,dA(w)
= \frac{1}{1-|z|^2} \int_D \frac{|(S_z 1)(\phi_z(w))|\,|{\phi_z}'(w)|}{\sqrt{1 -
|w|^2}}\,dA(w).
\]
In the last integral, make the substitution
$w = \phi_z(\lambda)$ and use the identities
\begin{gather}
\phi_z(\phi_z(\lambda)) = \lambda \notag \\[9pt]
{\phi_z}'(\phi_z(\lambda))
= \frac{1}{{\phi_z}'(\lambda)} = \frac{(1 -\bar{z}\lambda)^2}{|z|^2 - 1} \notag\\[9pt]
1 - |\phi_z(\lambda)|^2
= \frac{(1 - |z|^2)(1 - |\lambda|^2)}{|1 - \bar{z}\lambda|^2} \label{78} \\[9pt]
dA(w) = |{\phi_z}'(\lambda)|^2\,dA(\lambda) = \frac{(1-|z|^2)^2}{|1
-\bar{z}\lambda|^4}\,dA(\lambda) \notag
\end{gather}
to obtain
\[
\int_D \frac{|(SK_z)(w)|}{\sqrt{1 - |w|^2}}\,dA(w)
= \frac{1}{\sqrt{1-|z|^2}}
\int_D\frac{|(S_z 1)(\lambda)|}{|1-\bar{z}\lambda|\sqrt{1-
|\lambda|^2}}\,dA(\lambda).
\]
H\"older's inequality now gives
\[
\int_D \frac{|(SK_z)(w)|}{\sqrt{1 - |w|^2}}\,dA(w)
\le \frac{\|S_z 1\|_6}{\sqrt{1-|z|^2}} \left(\int_D
\frac{dA(\lambda)}{|1-\bar{z}\lambda|^{6/5} (1-|\lambda|^2)^{3/5}}\right)^{5/6}.
\]
{}Lemma 4 of Axler's paper \bibref{Axl} asserts that there exists a constant $c <
\infty$ independent of~$z$ (and of course also independent of~$S$) such that the last
integral on the right is bounded by~$c$.  This completes the proof of~\eqref{50}.

To prove \eqref{51}, replace $S$ with $S^*$ in \eqref{50}, interchange the roles of
$w$ and $z$ in
\eqref{50} and then use the equation
\begin{equation} \label{53}
(S^* K_w)(z) = \ip{S^*K_w, K_z} = \ip{K_w, S K_z} = \overline{(SK_z)(w)}
\end{equation}%
to obtain the desired result.
\qed

The next lemma is the final tool needed for the proof of
Theorem~\ref{37}.

\begin{lemma} \label{62}
If\/ $S$ is a finite sum of operators of the form\/ $T_{u_1} \dots T_{u_n}$, where
each\/ $u_j \in L^\infty(D,dA)$, then
\[
\sup_{z \in D} \|S_z1\|_p < \infty
\]
for every\/ $p \in (1,\infty)$.
\end{lemma}

\proof
By the triangle inequality, we can assume without loss of
generality that $S = T_{u_1} \dots T_{u_n}$, where each $u_j \in L^\infty(D,dA)$. 
Using \eqref{26} and the definition of $S_z$, we have
\begin{equation} \label{63}
S_z = T_{u_1 \circ \phi_z} \dots T_{u_n \circ \phi_z}.
\end{equation}%

Fix $p \in (1,\infty)$.  Then there is a constant $c < \infty$ such that
\mbox{$\|Pv\|_p \le c \|v\|_p$} for every $v \in L^p(D, dA)$ (see \bibref{Ax2},
Theorem~1.10).  Thus
\mbox{$\|T_u f\|_p \le c \|u\|_\infty \|f\|_p$} for all $u \in L^\infty(D,dA), f \in
L^2_a$.  Because
$\|u \circ \phi_z\|_\infty$ is independent of $z$, this implies that
$\|T_{u \circ \phi_z} f\|_p \le c \|u\|_\infty \|f\|_p$ for all
$u \in L^\infty(D,dA)$, all $f \in L^2_a$, and all $z \in D$.  This, along with
\eqref{63}, shows that
$\|S_z1\|_p$ is bounded independent of~$z$, as desired.
\qed

\section{The Proof}

Now we are ready to prove Theorem~\ref{37}.  We will show that
\[
\text{(a) $\Rightarrow$ (b) $\Rightarrow$ (c) $\Rightarrow$ (d)
$\Rightarrow$ (e)  $\Rightarrow$ (f) $\Rightarrow$ (a).}
\]
Several of these implications are easy or trivial.  The difficult parts of the proof
are the implications (c) $\Rightarrow$ (d) and (f) $\Rightarrow$ (a).

\bigskip

\textsc{Proof of Theorem~\ref{37}:} \hspace{.4em minus .2em}
First suppose that (a) holds, so $S$ is compact.  Because $k_z \to 0$ weakly in
$L^2_a$ as $z \to \partial D$, this implies that $\|Sk_z\|_2 \to 0$ as
$z \to \partial D$, completing the proof that
(a) implies~(b).

Now suppose that (b) holds, so $\|Sk_z\|_2 \to 0$ as
$z \to \partial D$.  Then
\[
|\tl{S}(z)| = |\ip{Sk_z, k_z}| \le \|Sk_z\|_2.
\]
Hence $\tl{S}(z) \to 0$ as $z \to \partial D$, completing the proof that (b)
implies~(c).

Now suppose that (c) holds, so $\tl{S}(z) \to 0$ as $z \to \partial D$.  To
prove~(d), it suffices to show that $\ip{S_z 1, w^n} \to 0$ as $z \to \partial D$ for
every nonnegative integer~$n$.  So fix a nonnegative integer $n$.

For $z, \lambda \in D$, we have
\begin{align*}
\tl{S}(\phi_z(\lambda)) &= \wt{S_z}(\lambda) \\
&= (1-|\lambda|^2)^2
\sum_{j,m=0}^\infty (j+1)(m+1) \ip{S_z w^j, w^m} \bar{\lambda}^j \lambda^m,
\end{align*}
where the first equality comes from Lemma~\ref{77} and the second equality comes
from~\eqref{70}.  Now fix $r \in (0,1)$, multiply both sides of the last equation by
$\bar{\lambda}^n / (1-|\lambda|^2)^2$, and then integrate over $rD$ to obtain
\begin{align*}
\int_{rD} \frac{\tl{S}(\phi_z(\lambda))
\bar{\lambda}^n}{(1-|\lambda|^2)^2}\,dA(\lambda) &= \sum_{j,m=0}^\infty (j+1)(m+1)
\ip{S_z w^j, w^m}
\int_{rD}\bar{\lambda}^{j+n} \lambda^m\,dA(\lambda) \\
&= \sum_{j=0}^\infty (j+1) \ip{S_zw^j, w^{j+n}} r^{2j+2n+2} \\
&= r^{2n+2} \left( \ip{S_z1,w^n} + \sum_{j=1}^\infty (j+1) \ip{S_zw^j, w^{j+n}}
r^{2j} \right).
\end{align*}
For each $\lambda \in D$, we know that $\phi_z(\lambda) \to \partial D$ as $z \to
\partial D$ (this follows from~\eqref{78}).  Thus by our hypothesis (c),
$\tl{S}(\phi_z(\lambda)) \to 0$ as $z \to \partial D$ for each $\lambda \in D$.  Hence
for each fixed $r \in (0,1)$, the left side of the equality above has limit $0$ as $z
\to \partial D$ (note that the integrand is bounded by $\|S\|/(1-r^2)^2$,
independently of $\lambda$ and~$z$, justifying our passage to of the limit).  Dividing
the right side of the last equality above by $r^{2n+2}$ (we are thinking of $r$ as
fixed), we conclude that
\begin{equation} \label{80}
\ip{S_z1,w^n} + \sum_{j=1}^\infty (j+1) \ip{S_zw^j, w^{j+n}} r^{2j} \to 0
\text{\ as\ } z \to \partial D
\end{equation}%
for each $r \in (0,1)$.

Let's examine the infinite sum above.  Note that
\begin{align}
\bigl|\sum_{j=1}^\infty (j+1) \ip{S_zw^j, w^{j+n}} r^{2j}\bigr|
&\le \|S\| \sum_{j=1}^\infty r^{2j} \label{79} \\
&= \|S\| \frac{r^2}{1-r^2}. \notag
\end{align}
Thus given $\epsilon > 0$, we can choose $r \in (0,1)$ such that the left side of
\eqref{79} is less than $\epsilon$ for all $z \in D$.  Hence \eqref{80} implies that
\[
\limsup_{z \to \partial D} | \ip{S_z1, w^n} | \le \epsilon.
\]
Because $\epsilon$ is an arbitrary positive number, the inequality above implies that
$\ip{S_z1, w^n} \to 0$ as $z \to \partial D$, as desired, completing the proof that
(c) implies~(d).

Now suppose that (d) holds, so $S_z 1 \to 0$ weakly in $L^2_a$ as
$z \to \partial D$.  For every $z \in D$ and every $r \in (0,1)$ we have
\begin{align}
{\|S_z 1\|_2}^2
&= \int_{D\setminus r\bar{D}} |(S_z 1)(w)|^2\,dA(w)
+ \int_{r\bar{D}} |(S_z 1)(w)|^2\,dA(w) \notag \\
& \le (1-r^2)^{1/2} {\|S_z 1\|_4}^2
+ \int_{r\bar{D}} |(S_z 1)(w)|^2\,dA(w), \label{76}
\end{align}
where the inequality comes from writing $|(S_z 1)(w)|^2 = 1 |(S_z 1)(w)|^2$ and then
using the Cauchy-Schwarz inequality in the first integral above.  By Lemma~\ref{62},
$\|S_z 1\|_4$ is bounded independent of~$z$.  Thus given $\epsilon > 0$, we can choose
$r \in (0,1)$ such that the first term in \eqref{76} is less than $\epsilon/2$ for all
$z \in D$.  Having chosen $r \in (0,1)$, the second term in \eqref{76} will also be
less than
$\epsilon/2$ for all $z$ sufficiently close to $\partial D$ (because a
sequence converging weakly to $0$ in $L^2_a$ converges uniformly to $0$ on each
compact subset of~$D$).  Thus for $z$ sufficiently close to $\partial D$
we have ${\|S_z 1\|_2}^2 < \epsilon$, completing the proof that (d) implies (e).

Now suppose that (e) holds, so
\begin{equation} \label{68}
\|S_z1\|_2 \to 0 \text{\ as\ } z \to \partial D.
\end{equation}%
To prove (f), fix $p \in (1, \infty)$.  If $1 < p \le 2$, then clearly \eqref{68}
implies that $\|S_z1\|_p \to 0$ as $z \to \partial D$.  To consider
the remaining case, suppose now that $2 < p < \infty$.  Then
\begin{equation} \label{61}
\|S_z1\|_p \le {\|S_z1\|_2}^{1/p} {\|S_z1\|_{2p-2}}^{(p-1)/p},
\end{equation}%
as can be seen by writing $|S_z1|^p = |S_z1| |S_z1|^{p-1}$ and then using
the Cauchy-Schwarz inequality in the integral defining $\|S_z1\|_p$.  By our
hypothesis~(e), the first term on the right side of \eqref{61} has limit $0$ as
$z \to \partial D$.  The second term on the right side of \eqref{61} is bounded
independent of $z$, by Lemma~\ref{62}.  Thus the left side of \eqref{61} has limit $0$
as $z \to \partial D$, completing the proof that (e) implies~(f).

Now suppose that (f) holds, so $\|S_z1\|_6 \to 0$ as
$z \to \partial D$.  For $f \in L^2_a$ and $w \in D$, we have
\begin{align}
(Sf)(w) &= \ip{Sf, K_w} \notag \\
&= \ip{f, S^*K_w} \notag \\
&= \int_D f(z) \overline{(S^*K_w)(z)}\,dA(z) \notag \\
&= \int_D f(z) (SK_z)(w)\,dA(z), \label{55}
\end{align}
where the last equation follows from \eqref{53}.

For $0 < r < 1$, define an operator $S_{[r]}$ on $L^2_a$ by
\begin{equation} \label{54}
(S_{[r]}f)(w) = \int_{rD} f(z) (SK_z)(w)\,dA(z).
\end{equation}%
In other words, $S_{[r]}$ is the integral operator with kernel
$(SK_z)(w) \chi_{rD}(z)$.  For each $r \in (0,1)$ we have
\begin{align*}
\int_D \int_D |(SK_z)(w) \chi_{rD}(z)|^2\,dA(w)\,dA(z)
&= \int_{rD} \int_D |SK_z(w)|^2\,dA(w)\,dA(z) \\
&= \int_{rD} {\|SK_z\|_2}^2\,dA(z) \\
&\le \|S\|^2 \int_{rD} {\|K_z\|_2}^2\,dA(z) \\
&< \infty.
\end{align*}
Thus $S_{[r]}$ is a Hilbert-Schmidt operator and in particular is compact.  Hence to
prove that (f) implies (a), we only need show that
\mbox{$\|S - S_{[r]}\| \to 0$} as $r \to 1^-$.

If $r \in (0,1)$, then $S - S_{[r]}$ is the integral operator with kernel
\[
(SK_z)(w) \chi_{D \setminus rD}(z),
\]
as can be seen from \eqref{55} and~\eqref{54}.  The Schur test (see page~126
of~\bibref{BHS}) implies that if $u$ is a positive measurable function on $D$ and
$c_1,  c_2$ are constants such that
\begin{equation} \label{56}
\int_D |(SK_z)(w) \chi_{D \setminus rD}(z)| u(w)\,dA(w) \le c_1 u(z)
\end{equation}%
for all $z \in D$ and
\begin{equation} \label{57}
\int_D |(SK_z)(w) \chi_{D \setminus rD}(z)| u(z)\,dA(z) \le c_2 u(w)
\end{equation}%
for all $w \in D$, then
\begin{equation} \label{64}
\|S - S_{[r]}\| \le \sqrt{c_1 c_2}.
\end{equation}%
Note that the left side of \eqref{56} equals $0$ for $0 < |z| < r$.  Taking $u(\lambda)
= 1/\sqrt{1-|\lambda|^2}$, we see from Lemma~\ref{58} that \eqref{56} and \eqref{57}
are satisfied with
$c_1 = c\,\sup\{\|S_z1\|_6: r \le |z| < 1\}$ and
$c_2 = c\,\sup\{\|S^*_w1\|_6: w \in D\}$, where $c$ is the constant from
Lemma~\ref{58}.  Our hypothesis (f) implies that $c_1 \to 0$ as
$r \to 1^-$, and Lemma~\ref{62} (with $S^*$ replacing $S$) shows that
$c_2 < \infty$.  Thus from \eqref{64} we conclude that
$\|S - S_{[r]}\| \to 0$ as $r \to 1^-$, completing the proof that (f)
implies~(a).
\qed

\section*{References}

\begin{list}{\arabic{referencec}.\hfill}{\usecounter{referencec} 
\setlength{\topsep}{12pt plus 3pt minus 3pt}
\setlength{\partopsep}{0pt}
\setlength{\labelwidth}{1.5\parindent}
\setlength{\labelsep}{0pt}
\setlength{\leftmargin}{1.5\parindent}
\setlength{\parsep}{1ex} }

\raggedright

\item \label{Ax2}
Sheldon Axler,
Bergman spaces and their operators,
\textsl{Surveys of Some Recent Results in Operator Theory}, vol.~1, edited by John
B.~Conway and Bernard B.~Morrel,
Pitman Research Notes in Mathematics, 1988, 1--50.

\item \label{Axl}
Sheldon Axler,
The Bergman space, the Bloch space, and commutators of multiplication operators,
\textsl{Duke Math. J.} 53 (1986), 315--332.

\item \label{AxC}
Sheldon Axler and \u{Z}eljko \u{C}u\u{c}kovi\'c,
Commuting Toeplitz operators with harmonic symbols,
\textsl{Integral Equations Operator Theory} 14 (1991), 1--12.

\item \label{AxZ}
Sheldon Axler and Dechao Zheng,
The Berezin transform on the Toeplitz algebra,
\textsl{Studia Math.} 127 (1998), 113--136.

\item \label{BHS}
Arlen Brown, P. R. Halmos, and A. L. Shields,
Ces\`aro operators,
\textsl{Acta Sci. Math. (Szeged)} 26 (1965), 125--137.

\item \label{KZh}
Boris Korenblum and Kehe Zhu,
An application of Tauberian theorems to Toeplitz operators,
\textsl{J. Operator Theory} 33 (1995), 353--361.

\item \label{NoR}
Eric Nordgren and Peter Rosenthal,
Boundary values of Berezin symbols,
\textsl{Operator Theory: Advances and Applications} 73 (1994), 362-368.

\item \label{St2}
Karel Stroethoff,
Compact Hankel operators on the Bergman space,
\textsl{Illinois J. Math.} 34 (1990), 159--174.

\item \label{Str}
Karel Stroethoff,
Compact Toeplitz operators on Bergman spaces,
\textsl{Math. Proc. Cambridge Philos. Soc.} 124 (1998), 151--160.

\item \label{St3}
Karel Stroethoff,
The Berezin transform and operators on spaces of analytic functions,
\textsl{Banach Center Publ.} 38 (1997), 361--380.

\item \label{StZ}
Karel Stroethoff and Dechao Zheng,
Toeplitz and Hankel operators on Bergman spaces,
\textsl{Trans. Amer. Math. Soc.} 329 (1992), 773--794.

\item \label{Zhe}
Dechao Zheng,
Toeplitz and Hankel operators,
\textsl{Integral Equations Operator Theory} 12 (1989), 280--299.

\item \label{Zhu}
Kehe Zhu,
Positive Toeplitz operators on weighted Bergman spaces of bounded symmetric domains,
\textsl{J. Operator Theory} 20 (1988), 329--357.

\end{list}

\bigskip

\bigskip

\pagebreak[3]
\noindent
{\scshape
Sheldon Axler} \hfill \textsl{e-mail}: \texttt{axler@sfsu.edu}\\
{\scshape
Department of Mathematics}
\hfill \textsl{www}:
\texttt{http://math.sfsu.edu/axler}\\
{\scshape
San Francisco State University \\
San Francisco, CA 94132 USA

\bigskip

\noindent
Dechao Zheng} \hfill \textsl{e-mail}: \texttt{zheng@math.vanderbilt.edu}\\
{\scshape
Department of Mathematics \\
Vanderbilt University \\
Nashville, TN 37240 USA}

\end{document}